\documentclass[12pt,oneside]{amsart}
\usepackage{cite}
\usepackage{latexsym}
\usepackage{layout}
\usepackage{amsfonts}
\usepackage{amssymb}
\usepackage{epsfig}

\newtheorem{thm}{Theorem}[section]

\newtheorem{deff}[thm]{Definition}
\newtheorem{lem}[thm]{Lemma}
\newtheorem{cor}[thm]{Corollary}

\newtheorem{prp}[thm]{Proposition}

\begin{document}
\title{On the quadratic formula modulo $N$}
\author{Steve 
Wright} 
%\date{\today}
\address{Steve Wright, Department of Mathematics and Statistics,
Oakland University, Rochester, MI 48309. Email: wright@oakland.edu}

\begin{abstract}
Let $a, b, c,$ and $n$ be integers, with $a$ nonzero and $n$ at least two. Necessary and sufficient conditions on these parameters are derived which guarantee that all solutions of the congruence
\[
ax^2+bx+c \equiv 0\ \textrm{mod}\ n
\]
are given precisely by the solutions of
\[
2ax\equiv -b+s \ \textrm{mod}\ n,
\]
where $s$ varies over all solutions of 
\[
x^2\equiv b^2-4ac \ \textrm{mod}\ n.
\]
Corollaries of this result are deduced for prime-power moduli and some illustrative examples are also presented.
 \end{abstract}
\maketitle
\markboth{}{}
\noindent \textit{keywords}: \textrm{quadratic formula, quadratic residue, quadratic non-residue, congruence modulo $n$, quadratic congruence}

\section{Introduction}

Let $a, b, c,$ and $n$ be fixed integers, with $a$ nonzero and $n$ at least two. In Section IV of the \emph{Disquisitiones Arithmeticae}, Gauss presented a complete and elegant solution of the quadratic congruence
\begin{equation*}
ax^2+bx+c \equiv 0\ \textrm{mod}\ n. \tag{1.1}
\]
By completing the square in $ax^2+bx+c $, it is easy to see that if $d=b^2-4ac$ is the discriminant of the quadratic, then the solutions of (1.1) are obtained as the solutions of
\begin{equation*}
2ax\equiv -b+s \ \textrm{mod}\ 4an, \tag{1.2}
\]
where $s$ varies over all solutions of
\begin{equation*}
x^2\equiv d \ \textrm{mod}\ 4an \tag{1.3}
\]
that are pairwise incongruent mod $2an$. Thus the solution of (1.1) is reduced to the solution of the ``pure" quadratic congruence (1.3) (Gauss' terminology), and it is the solution of this latter congruence, which we will call the (modular) \emph{square-root problem}, that Gauss devotes his attention to in the \emph{Disquisitiones}. Needless to say, Gauss' work here is a milestone of number theory, and has been a guide and inspiration to the subject ever since.

By setting $s=\sqrt{b^2-4ac}$, one may write (1.2) as
\[
2ax\equiv -b+\sqrt{b^2-4ac} \ \textrm{mod}\ 4an,
\]
which is reminiscent of the quadratic formula for quadratic equations from secondary-school algebra. If one wants an \emph{exact} analog of the quadratic formula, one would require that all solutions of (1.1) are determined from the equation
\begin{equation*}
2ax\equiv -b+s \ \textrm{mod}\ n, \tag{1.4}
\]
where $s$ varies over all solutions of 
\begin{equation*}
x^2\equiv d \ \textrm{mod}\ n. \tag{1.5}
\]
One would then want to find all solutions of (1.4) by simply ``dividing" by $2a$, i.e., multiplying by an inverse of $2a$ mod $n$ to obtain
\begin{equation*}
x\equiv \frac{ -b+\sqrt{b^2-4ac}}{2a}\ \textrm{mod}\ n. \tag{1.6}
\]
Since this requires the existence of the indicated inverse, a necessary condition for the solutions of (1.1) to be given by (1.6) is to have $2a$ and $n$ relatively prime, and a simple argument shows that this is also sufficient. We will refer to the solution of (1.1) that is given by (1.2) and (1.3) as the \emph{general form} of the quadratic formula, and we will call the solution of (1.1) given by (1.4)-(1.6) the \emph{exact form}.

When it can be applied, the exact form of the quadratic formula is obviously a more efficient way to solve (1.1) than the general form; the disadvantage is that is requires the rather restrictive condition of the relative primality of $2a$ and $n$. A question which thus naturally occurs asks if (1.4) and (1.5) (but not necessarily (1.6)) can be used to solve (1.1) \emph{without} this condition, and if so, to determine precisely for which moduli $n$ this can be done. We will say that the \emph{intermediate form} of the quadratic formula (IQF) is valid if (1.4) and (1.5) completely solve (1.1).  The purpose of this paper is to show that IQF can indeed hold when $2a$ and $n$ are not relatively prime and to characterize precisely the moduli for which it is valid. The answer is given by Theorem 6.1 in Section 6 (see also Definition 3.4 in Section 3) and is, at least to us, surprisingly subtle. Moreover, our methods are entirely elementary; indeed, everything required for our analysis (and much more!) is already contained in the \emph{Disquisitiones}.

We will now briefly describe the contents of the paper. In Section 2, we introduce notation and terminology (most of which is quite standard) that will be used throughout the sequel and state the results on which the rest of our work depends, the most essential of which is Gauss' solution of the square-root problem. The analysis of IQF begins in Section 3, where it is reduced to two statements relating the set of all solutions of (1.1) to the set of all solutions of an associated congruence. Three results required for the study of these solution sets are also established. Necessary and sufficient conditions for the reduction in Section 3 to be valid are derived in Sections 4 and 5. Section 6 contains the main result, Theorem 6.1, which is an immediate consequence of the work of the previous three sections. Two corollaries for prime-power moduli are deduced from it, and some illustrative examples are also presented.
\section{Preliminaries}

We begin with some notation and terminology that will be used systematically throughout the remainder of this paper. Let $\mathbf{Z}$ denote the set of integers, and $\mathbf{Z}^{+}$ the set of positive integers. The symbol $\emptyset$ will denote the empty set. If $p$ is a prime number and $z$ is an integer, we will let $\mu_p(z)$ denote the multiplicity of $p$ in $z$, and take $\mu_p(z)=0$ if $p$ is not a factor of $z$. If $a$ and $b$ are integers, then $(a, b)$ will denote the greatest common divisor of $a$ and $b$. For $a, b, c\in \mathbf{Z}$, we set 
\[
q(x)=ax^2+bx+c,
\]
and let $d=b^2-4ac$ denote the discriminant of $q(x)$.

If $n$ is a positive integer, we will say that an integer $a$ is a quadratic residue or non-residue of $n$ if the equation $x^2\equiv a\ \textrm{mod}\ n$ either does, or does not, have a solution $x$ in $\mathbf{Z}$. The set of quadratic residues of $n$ will be denoted by $Q(n)$. The following proposition will prove quite useful to us, and is a simple consequence of the difference-of-two-squares factorization identity and the Chinese remainder theorem. We note here that in all of what follows, \emph{a solution to a modular congruence will always mean a solution that is nonnegative and minimal with respect to the relevant modulus}, i.e., if $n$ is the modulus and $\sigma$ is a solution, then $0\leq \sigma < n$.
\begin{prp}
If $k, l \in \mathbf{Z}^{+}, (k, l)=1,$ and $a \in \mathbf{Z}$, then $\sigma$ is a solution of  $x^2\equiv a\ \textnormal{mod}\ (kl)$ if and only if there exist solutions $\kappa$ and $\lambda$ of $x^2\equiv a\ \textnormal{mod}\ k$ and $x^2\equiv a\ \textnormal{mod}\ l$, respectively, such that $\sigma \equiv \kappa \ \textnormal{mod}\ k$ and $\sigma \equiv \lambda \ \textnormal{mod}\ l$.
\end{prp}

\begin{cor}
If $k, l \in \mathbf{Z}^{+}$ and $(k, l)=1$, then
\[
Q(kl)=Q(k)\cap Q(l).
\]
\end{cor}

Our study of the intermediate form of the quadratic formula will make essential use of Gauss' beautiful solution of the square-root problem as set forth in \emph{Disquisitiones Arithmeticae}. We will now describe this solution in detail.

Let $p$ be a fixed prime, $k \in \mathbf{Z}^{+}, u \in \mathbf{Z}$. We suppose first that $u \in Q(p^k)$ and consider solutions $\sigma$ of the congruence $x^2\equiv u\ \textrm{mod}\ p^k$. In [1, article 104], we find these solutions determined as follows:

I. Suppose first that $u$ is not divisible by $p$. If $p=2$ and $k=1$ then $\sigma=1$. If $p$ is odd or $p=2=k$ then $\sigma$ has exactly two values $\pm \sigma_0$. Finally, if $p=2$ and $k>2$ then $\sigma$ has exactly four values  $\pm \sigma_0$ and $\pm \sigma_0+2^{k-1}$. 

II. If $u$ is divisible by $p$ but not by $p^k$, let $2\mu=\mu_p(u)$ (which necessarily must be even when $u \in Q(p^k)$) and  let $u=u_1p^{2\mu}$. Then $\sigma$ is given by the formula
\[
\sigma^{\prime}p^{\mu}+ip^{k-\mu},\ i\in \{0, 1,\dots,p^{\mu}-1\},\]
where $\sigma^{\prime}$ varies over all solutions, determined according to I, of the congruence
\[
x^2 \equiv u_1\ \textrm{mod}\ p^{k-2\mu}.\]

III. If $u$ is divisible by $p^k$, and if we set $k=2m$ or $k=2m-1$, depending on whether $k$ is even or odd, then $\sigma$ is given by the formula
\[
ip^m,\ i\in \{0,\dots,p^{k-m}-1\}.\]

If $v$ is now an arbitrary modulus greater than 1 and $u\in Q(v)$, then the solutions $\sigma$ of $x^2\equiv u\ \textrm{mod}\ v$ are given precisely via the prime factorization $p_1^{\alpha_1}\cdots p_t^{\alpha_t}$ of $v$ and Proposition 2.1 by the recipe
\[
\sigma \equiv\ u_i \ \textrm{mod}\ p_i^{\alpha_i},\]
where $u_i$ is any solution, determined according to I, II, or III, of
\[
x^2\equiv u\ \textrm{mod}\ p^{\alpha_i}, i=1,\dots, t.
\]
We will refer to all of this as \emph{Gauss' solution of the square-root problem}.
\section{Analysis of IQF: the Initial Reduction}

In this section we determine a condition equivalent to the validity of IQF that relates the solution set of $q(x) \equiv\ 0\ \textrm{mod}\ n$ to the solution set of an auxiliary congruence. We also establish some lemmas that will be used to study this relationship more closely.

Let $a, b , c,$ and $n$ be fixed integers with $n>1$ and $a$ nonzero. Let $d=b^2-4ac$ and $q(x)=ax^2+bx+c$. In all of what follows, the phrase ``IQF is true" will mean that IQF is true for the congruence $q(x)\equiv\ 0\ \textrm{mod}\ n$. Completion of the square in $q(x)$ shows that
\begin{quote}
(3.1) \hspace {0.5 cm} IQF is true if and only if for all $x\in \mathbf{Z}, 4aq(x)\equiv\ 0\ \textrm{mod}\ n$ if and \hspace{3cm} 
\end{quote}
\hspace{2.9cm} only if $q(x) \equiv\ 0\ \textrm{mod}\ n$.

Now let $r=(a, n), a_1=a/r,$ and $k=$ multiplicity of 2 in $n/r$. Then $n=2^krm$, where $m$ is odd and $(a_1, 2^km)=1$. In particular, $a_1$ is odd if $k>0$ and $(2a, n)>1$ if and only if either $r>1$ or $k>0$.
\begin{lem}
Let z $\in \textbf{Z}$. If $k=0, 1,$ or $2$, then $4az$ is divisible by n if and only if z is divisible by m. If $k\geq 3$, then $4az$ is divisible by n if and only if z is divisible by $2^{k-2}m$.
\end{lem}

\emph{Proof}. If $z \in \textbf{Z}$ then $4az$ is divisible by $n$ if and only if $4a_1z$ is divisible by $m$ if $k=0$, $2a_1z$ is divisible by $m$ if $k=1$, $a_1z$ is divisible by $m$ if $k=2$, or $a_1z$ is divisible by $2^{k-2}m$ if $k\geq 3$. Since $m$ is odd, $(a_1, m)=1$, and $a_1$ is odd for $k>0$, it follows that $(m, a_1)=(m, 2a_1)=(m, 4a_1)=1$, and $ (2^{k-2}m, a_1)=1$ if $k\geq 3$. The conclusions of the lemma are now simple consequences of all of this. $\hspace{10.7cm} \textrm{QED}$

Let
\[
Q=\left\{\begin{array}{cc}\{x \in \textbf{Z}: q(x) \equiv\ 0\ \textrm{mod}\ m\},\ \textrm{if $k=0, 1,$ or $2$,}\\
\{x \in \textbf{Z}: q(x) \equiv\ 0\ \textrm{mod}\ 2^{k-2}m\},\ \textrm{if $k\geq$ 3,}\\\end{array}\right.
\]
\vspace{0.3cm}
\[
T=\{x \in \textbf{Z}: q(x) \equiv\ 0\ \textrm{mod}\ n\}.\]
It is now an immediate consequence of (3.1) and Lemma 3.1 that
\begin{center}
IQF is true if and only if $Q=T$.
\end{center}
In light of this observation and the fact that $T\subseteq Q$, IQF will thus be valid if and only if either
\begin{equation*}
Q=\emptyset, \tag{3.2}\]
or
\begin{equation*}
\emptyset\not= Q=T. \tag{3.3}\]
The derivation of necessary and sufficient conditions which guarantee the validity of (3.2) and (3.3) will be carried out in Sections 4 and 5, respectively.

The following lemma will play a pivotal role in our analysis of (3.3) in Section 5. In order to state it, we first let $u, v\in \textbf{Z}^{+}$, with $q(x)$ and $d$ as specified at the beginning of this section. If $\mathcal{S}_0$ (respectively, $\mathcal{S}_1$) denotes the set of all solutions of $x^2\equiv\ d\ \textrm{mod}\ 4auv$ (respectively, $x^2\equiv\ d\ \textrm{mod}\ 4av$) that are pairwise incongruent mod $2auv$ (respectively, mod $2av$), then we set $\Sigma_i=\{\sigma \in \mathcal{S}_i: \sigma \equiv b\ \textrm{mod}\ 2a\}, i=0, 1$. We note that $\Sigma_0$ (respectively, $\Sigma_1$) is uniquely determined up to congruence mod $2auv $(respectively, mod $2av$).

\begin{lem}
If $u, v\in \textbf{Z}^{+}$, $q(x)$ and d are as specified at the beginning of this section, $\Sigma_0$ and $\Sigma_1$ are as defined above,
\[
Q_0=\{x \in \textbf{Z}: q(x)\equiv 0\ \textnormal{mod}\ uv\},\]
\[
Q_1=\{x \in \textbf{Z}: q(x)\equiv 0\ \textnormal{mod}\ v\},\ and\]
\begin{quote}
 $S_i=$the set of all elements of $Q_i$ minimal and nonnegative with respect to the appropriate modulus, $i=0, 1$,
\end{quote}
then the following statements are equivalent:

$(a)\ \emptyset \not= Q_1=Q_0$;

$(b)\ \emptyset \not= S_1$ and $S_0=\big\{s+jv: s \in S_1, j\in \{0,\dots,u-1\}\big\}$;

$(c)\ \emptyset \not= \Sigma_1$ and for each $\sigma\in \Sigma_1$ and $j\in \{0,\dots,u-1\}$, there exista $\sigma^{\prime}\in \Sigma_0$ such that
\[
\sigma^{\prime} \equiv\ \sigma+2avj\ \textnormal{mod}\ 2auv.\]
Furthermore, if $(2a, uv)=1$ and if
\[
\Sigma_0^{\prime}=\ \textrm{the set of all solutions of}\ x^2\equiv\ d \ \textnormal{mod}\ uv,\] 
\[
\Sigma_1^{\prime}=\ \textrm{the set of all solutions of}\ x^2\equiv\ d \ \textnormal{mod}\ v,\] 
then $(a), (b),$ and $(c)$ are equivalent to

$(d)\ \emptyset \not=\Sigma_1^{\prime} $ and $\Sigma_0^{\prime}=\big\{\sigma+jv: \sigma \in \Sigma_1^{\prime},\ j\in \{0,\dots,u-1\}\big\}$.
\end{lem}

\emph{Proof}. $(a)\Rightarrow (b)$ Let $S_2$ denote the set on the right-hand side of the equation in $(b)$. Then $S_2\subseteq [0, uv)$. We have by $(a)$ that
\[
S_2\subseteq S_1+v\textbf{Z} =Q_1=Q_0.\]
But $ S_1+v\textbf{Z} \subseteq S_2+uv\textbf{Z}$, hence
\[
Q_0=Q_1\subseteq S_2+uv\textbf{Z}\subseteq Q_0,\]
i.e., $Q_0=S_2+uv\textbf{Z}$, and $(b)$ is an immediate consequence of this.

$(b)\Rightarrow (a)$ Clearly $Q_1\not= \emptyset$ and $Q_0\subseteq Q_1$. Hence from $(b)$, we obtain
\[
Q_1=S_1+v\textbf{Z} \subseteq S_2+uv\textbf{Z}=S_0+uv\textbf{Z}=Q_0.\]

$(b)\Rightarrow (c)$ By the general form of the quadratic formula, the elements of $S_0$ (respectively, $S_1$) consist precisely of the nonnegative minimal residues mod $uv$ (respectively, mod $v$) of
\[
\frac{\sigma-b}{2a},\ \sigma \in \Sigma_0\ (\textrm{respectively}, \sigma \in \Sigma_1)\]
(here we mean ordinary division and not multiplication by an inverse relative to the modulus). We evidently have $\Sigma_1\not= \emptyset$, so let $\sigma \in \Sigma_1$ and $j \in \{0,\dots,u-1\}$. Then there exists $s \in S_1$ such that
\[
 \frac{\sigma-b}{2a}\equiv s\ \textnormal{mod}\ v,\]
 hence one may find $j^{\prime}\in \{0,\dots,u-1\}$ such that
 \begin{equation*}
  \frac{\sigma-b}{2a}+jv\equiv s+j^{\prime}v\ \textnormal{mod}\ uv. \tag{3.4}\]
Now from $(b), s+j^{\prime}v \in S_0$, and so there is a $\sigma^{\prime}\in \Sigma_0$ such that
 \begin{equation*}
 s+j^{\prime}v\equiv \frac{\sigma^{\prime}-b}{2a}\ \textnormal{mod}\ uv. \tag{3.5}\]
It now follows from (3.4) and (3.5) that
\[
\sigma^{\prime} \equiv\ \sigma+2avj\ \textnormal{mod}\ 2auv.\]

$(c)\Rightarrow (b)$ Clearly $S_1\not= \emptyset$. If $S_2$ is as it was before, then $S_0 \subseteq S_2$. In order to verify the reverse inclusion take $s+jv\in S_2$ and find $\sigma \in \Sigma_1, j^{\prime}\in \{0,\dots,u-1\}$ for which
\[
 s+jv\equiv \frac{\sigma-b}{2a}+j^{\prime}v\ \textnormal{mod}\ uv.\]
By $(c)$, there exists $\sigma^{\prime} \in \Sigma_0$ such that
\[
s+jv\equiv \frac{\sigma^{\prime}-b}{2a}\ \textnormal{mod}\ uv,\]
and so $s+jv \in Q_0$. Since $0\leq s+jv<uv$, it follows that it must also be in $S_0$.

Next, suppose that $(2a, uv)=1$. We will show that $(b)$ is equivalent to $(d)$. Since $(2a, uv)=1$, the exact form of the quadratic formula shows that there is a bijection between $S_i$ and $\Sigma_i^{\prime}, i=0, 1$. If $(b)$ is true, then $\Sigma_0^{\prime}$ hence has the same cardinality as the set $\Sigma_2^{\prime}$ on the right-hand side of the equation in $(d)$. But the inclusions
\[
\Sigma_0^{\prime} \subseteq \Sigma_1^{\prime}+v\textbf{Z} \subseteq \Sigma_2^{\prime}+uv\textbf{Z}\]
hold, and so $\Sigma_0^{\prime} \subseteq \Sigma_2^{\prime}$, since both sets are contained in $[0, uv)$. Since $\Sigma_0^{\prime}$ and $\Sigma_2^{\prime}$ have the same (finite) cardinality, they must hence be equal. An exchange of the roles of $S_i$ and $\Sigma_i^{\prime}, i=0, 1$ in this argument proves that $(b)$ is a consequence of $(d)$. $\hspace{5.5cm} \textrm{QED}$

The next two results will provide us with the tools we need to derive conditions which insure the validity of (3.2). The first gives  necessary and sufficient conditions for a quadratic congruence to have no solutions and the second is a quadratic residue calculation that will prove useful.
\begin{prp}
Let a, b, c, n, d and $q(x)$ be as specified at the beginning of this section. The congruence $q(x)\equiv\ 0\ \textnormal{mod}\ n$ has no solutions if and only if either

$(a)$ d is a quadratic non-residue of $4an$, or

$(b)$ d is a quadratic residue of $4an$ and there exists a prime factor p of $2a$ with the following properties: if $\alpha$ is the multiplicity of p in $4an$ and $\beta$ is the multiplicity of p in $2a$, then

$(i)\ 1<\beta<\alpha$;

$(ii)$ b is divisible by p and d is divisible by $p^2$;

$(iii)$ if d is not divisible by $p^{\alpha}, 2\mu$ is the multiplicity of p in d, $d=d_1p^{2\mu}$ and $\Sigma$ is the set of all solutions of $x^2 \equiv\ d_1\ \textnormal{mod}\ p^{\alpha-2\mu}$, then
\begin{equation*}
\sigma p^{\mu}+ip^{\alpha-\mu} \not \equiv b\ \textnormal{mod}\ p^{\beta},\ \forall\ \sigma \in \Sigma,\ \forall\ i\in \{0, 1,\dots,p^{\mu}-1\}; \tag{3.6}\]

$(iv)$ if d is divisible by $p^{\alpha}$ and s is chosen so that $\alpha=2s$ if $\alpha$ is even or $\alpha=2s-1$ if $\alpha$ is odd, then
\begin{equation*}
ip^s\not \equiv b\ \textnormal{mod}\ p^{\beta},\ \forall\ i\in \{0, 1,\dots,p^{\alpha-s}-1\}. \tag{3.7}\]
\end{prp}

\emph{Proof}. It follows from the general form of the quadratic formula that  $q(x)\equiv\ 0\ \textnormal{mod}\ n$ has no solutions if and only if either $(a)$ is true or

\vspace{0.3cm}
$(b)^{\prime}\ d\in Q(4an)$ and $y\not \equiv b$ mod $2a$ for every solution $y$ of $x^2 \equiv\ d\ \textnormal{mod}\ 4an$.

\vspace{0.3cm}

\noindent Prime factorization in concert with Proposition 2.1 shows that $(b)^{\prime}$ is equivalent to the statement

\vspace{0.3cm}
$(c)\ d\in Q(4an)$ and there is a prime factor $p$ of $2a$ with the following property: if $\alpha=\mu_p(4an), \beta=\mu_p(2a)$ then $y\not \equiv b\ \textnormal{mod}\ p^{\beta}$ for every solution $y$ of $x^2 \equiv\ d\ \textnormal{mod}\ p^{\alpha}$.
\vspace{0.3cm}

\noindent Thus it suffices to show that $(b)$ and $(c)$ are equivalent, and since $(b)$ obviously implies $(c)$ in light of Gauss' solution to the square-root problem, we need only establish the converse.

We hence assume that $(c)$ is true. Observe first that
\begin{equation*}
b^2\equiv\ d\ \textnormal{mod}\ 4a. \tag{3.8}\]
It follows that $\beta<\alpha$; otherwise $(c)$ would be false. Suppose that $p$ is odd. If $p$ does not divide $b$, it follows from (3.8) and [1, article 101] that there is a solution of $x^2 \equiv\ d\ \textnormal{mod}\ p^{\alpha}$ that is congruent to $b$ mod $p^{\beta}$, again contrary to $(c)$. Hence $b$ is divisible by $p$. Suppose that $p=2$ and $b$ is odd. Then $d$ is odd by (3.8), and so every solution $y$ of $x^2 \equiv\ d\ \textnormal{mod}\ 2^{\alpha}$ is also odd. If $\beta=1$ then $y\not \equiv b$ mod 2 for all such $y$, i.e., $y$ and $b$ have opposite parity, which they do not. Thus $\beta>1$. Now $\mu_2(4a)=1+\beta>2$, and, by (3.8), $d\in Q(2^{1+\beta})$. Hence $d\equiv\ 1$ mod 8 [1, article 103] and so by (3.8) and [1, articles 88 and 103], there is a solution of $x^2 \equiv\ d\ \textnormal{mod}\ 2^{\alpha}$  that is congruent to $b$ mod $2^{\beta}$, and hence $(c)$ is contradicted yet again. Thus $b$ is even if $p=2$. It follows that $p^2$ divides $d$, and so either $(iii)$ or $(iv)$ of $(b)$ must hold, each being simply a restatement of the conclusion of $(c)$ using the explicit solutions of $x^2 \equiv\ d\ \textnormal{mod}\ p^{\alpha}$ that result from Gauss' solution of the square-root problem. Suppose finally that $\beta=1$. Then we set $i=0$ in either (3.6) or (3.7) to conclude that either $\sigma p^{\mu}\not \equiv 0$ mod $p$ or $b\not \equiv 0$ mod $p$, neither of which can be true, since $p$ divides $b$ and $\mu>0$. Hence $\beta>1$. $\hspace{14.2cm} \textrm{QED}$
\begin{deff}
If $p$ is a prime number, $\alpha, \beta \in \textbf{Z}^{+}$, and $b, d\in \textbf{Z}$, then we will say that $(p^{\alpha}, p^{\beta})$ $\textnormal{forms a $(b, d)$-obstruction}$ if either condition $(b)(iii)$ or condition $(b)(iv)$ in Proposition $3.3$ holds for $p, \alpha, \beta, b,$ and $d$.
\end{deff}
\begin{lem}
Let a, n, m, and k be as specified at the beginning of this section.

$(a)$ If $k=0, 1,$ or $2$, then $Q(n)\cap Q(4a)=Q(4am)$;

$(b)$ If $k\geq 3$, then $Q(n)\cap Q(4a)=Q(2^kam)$.
\end{lem}
\emph{Proof}. If $r=(a, n), \rho=\mu_2(r), a_1=a/r$, and $\sigma=\mu_2(a_1)$, then we have the factorizations
\[
n=2^{k+\rho}r_1m,\ 4a=2^{\rho+\sigma+2}r_1s_1,\]
\[
2^{\varepsilon}am=2^{\varepsilon+\rho+\sigma}r_1s_1m,\ \varepsilon=2\ \textrm{or}\ k,\]
where $m, r_1,$ and $s_1$ are all odd and $(m, s_1)=1$. Using these facts, Corollary 2.2, and the prime factorizations of $n, 4a,$ and $2^{\varepsilon}am,$ we can find a subset $X$ of $\textbf{Z}$ such that
\begin{equation*}
Q(n)\cap Q(4a)=Q(2^{k+\rho})\cap Q(2^{\rho+\sigma+2})\cap X, \tag{3.9}\]
\begin{equation*}
Q(2^{\varepsilon}am)=Q(2^{\varepsilon+\rho+\sigma})\cap X,\ \varepsilon=2\ \textrm{or}\ k, \tag{3.10}\]

If $\varepsilon=2$ and $k\leq 2$, then
\[
Q(2^{k+\rho})\cap Q(2^{\rho+\sigma+2})=Q(2^{\rho+\sigma+2}),\]
and so from (3.9) and (3.10) it follows that 
\[
Q(n)\cap Q(4a)=Q(2^{\rho+\sigma+2})\cap X=Q(4am).\]

If $k\geq 3$ then $n/r$ is even. Since $(a_1, n/r)=1, a_1$ must be odd, and so $\sigma=0$. If $r$ is even, then $\rho>0$ and
\[
Q(2^{k+\rho})\cap Q(2^{\rho+2})=Q(2^{k+\rho}),\]
and hence from (3.9) and (3.10) we obtain
\[
Q(n)\cap Q(4a)=Q(2^{k+\rho})\cap X=Q(2^kam).\]
If $r$ is odd then $\rho=0$, and so
\[
Q(n)\cap Q(4a)=Q(2^k)\cap Q(2^2)\cap X=Q(2^k)\cap X=Q(2^kam).\]
$\hspace{15.4cm} \textrm{QED}$

We close this section by noting that if $b^2-4ac$ is a quadratic non-residue of $n$, it is also obviously a quadratic non-residue of $4an$. It follows that both the general form and the intermediate form of the quadratic formula will produce no solutions of $q(x)\equiv\ 0\ \textnormal{mod}\ n$, and so IQF is true in this situation. We record this observation as
\begin{lem}
Let d, n, and $q(x)$ be as specified at the beginning of this section. If $d \not \in Q(n)$ then $\textnormal{IQF}$ holds for $q(x)\equiv\ 0\ \textnormal{mod}\ n$.
\end{lem}
\section{$Q=\emptyset$}

With Proposition 3.3 and Lemma 3.5 in hand, it is now a simple matter to determine when $Q=\emptyset$.
\begin{lem}
Let a, b, d, k, m, n, and Q be as defined at the beginning of Section $3$. If $d\in Q(n)$ then $Q=\emptyset$ if and only if there is a prime factor p of $2a$ such that if $\beta=\mu_p(2a)$ and
\[
\alpha=\left\{\begin{array}{cc} \mu_p(4am),\ \textrm{if $k=0, 1,$ or $2$,}\\
\mu_p(2^kam),\ \textrm{if $k\geq$ $3$,}\\\end{array}\right.
\]
then $1<\beta<\alpha,$ b is divisible by p, d is divisible by $p^2$, and $(p^{\alpha}, p^{\beta})$ forms a $(b, d)$-obstruction.
\end{lem}

\emph{Proof}. By hypothesis, $d \in Q(n)$ and it is always the case that $d\in Q(4a)$, and so it follows from Lemma 3.5 that $d\in Q(4am)$ if $k=0, 1, $ or 2, and $d\in Q(2^kam)$ if $k\geq 3$. The conclusion of Lemma 4.1 is now a consequence of Proposition 3.3. $\hspace{6cm} \textrm{QED}$
\section{$\emptyset \not= Q=T$}

We begin this section by deriving necessary conditions for $\emptyset \not= Q=T$ to be valid. We will then prove that these conditions are also sufficient.

\begin{lem}
Let a, b, c, m, n, r, k, $q(x)$, Q, and T be as specified at the beginning of Section $3$, and let $\delta=(m, r)$. If $\emptyset \not= Q=T$, then
\begin{equation*}
b\ \textrm{and}\ c\ \textrm{are divisible by}\ r, \tag{5.1}\] 
and
\begin{equation*}
if\ k=0,\ then\ d/r^2\in Q(m)\ and\ either\ \delta=1\ or\ \delta\ is\ the\ product\ of\ distinct\ odd \tag{5.2}
\end{equation*}

\hspace{0.7cm} primes $p_1,\dots, p_t,$ each prime $p_i$ has even multiplicity $m_i$ in $m$, and $d/r^2$ is divisible 

\vspace{0.2cm}
\hspace{0.7cm} by $p_1^{m_1}\cdots p_t^{m_t}$;
\vspace{0.2cm}
\begin{equation*}
if\ k=1,\ then\ b/r\ is\ odd,\ either\ a/r\ or\ c/r\ is\ even, and\ d/r^2\ and\ \delta\ satisfy\ the \tag{5.3}\]
\hspace{1.1cm} conditions specified for them in $(5.2)$;
\begin{equation*}
if\ k\geq 2,\ then\ r/\delta\ and\ k\ are\ odd,\ k-1=\mu_2(d/r^2),\ d/(r^22^{k-1})\equiv\ 1\ \textnormal{mod}\ 8,\ and \tag{5.4}\]
\hspace{1.1cm} $d/r^2$ and $\delta$ satisfy the conditions specified for them in $(5.2)$.
\end{lem}

\emph{Proof}. We begin with the verification of (5.1). Let $x \in Q$, and deduce from the assumption $Q=T$ that for all $z\in \mathbf{Z}$,
\[
q(x+zm)\equiv\ 0\ \textrm{mod}\ 2^krm,\ \textrm{if}\ k=0, 1, \ \textrm{or}\ 2,\]
or
\[
q(x+2^{k-2}zm)\equiv\ 0\ \textrm{mod}\ 2^krm,\ \textrm{if}\ k\geq 3,\]
from whence it follows that for all $z\in \textbf{Z}$,
\[
\frac{q(x)}{m}+bz\equiv\  0\ \textrm{mod}\ r,\ \textrm{if}\ k=0, 1, \ \textrm{or}\ 2,\]
or
\[
\frac{q(x)}{2^{k-2}m}+bz\equiv\  0\ \textrm{mod}\ r,\ \textrm{if}\ k\geq 3.\]
Thus $r$ divides $b$ and $q(x)$ and so $r$ also divides $c=q(x)-ax^2-bx$.

If $a_1=a/r,\ b_1=b/r,\ c_1=c/r,\ q_1(x)=a_1x^2+b_1x+c_1,\ m_1=m/\delta,\ r_1=r/\delta$, and $\rho=\mu_2(r_1)$, then a simple argument using the facts that $m_1$ and $r_1/2^{\rho}$ are odd and $(m_1, r_1)=1$ confirms that if we set
\[
Q_0=\{x \in \textbf{Z}: q_1(x)\equiv\ 0\ \textrm{mod}\ 2^km\},\]

\[
Q_1=\left\{\begin{array}{cc}\{x \in \textbf{Z}: q_1(x) \equiv\ 0\ \textrm{mod}\ m_1\},\ \textrm{if $\rho\geq k-2$.}\\
\{x \in \textbf{Z}: q_1(x) \equiv\ 0\ \textrm{mod}\ 2^{k-\rho-2}m_1\},\ \textrm{if $0\leq \rho<k-2$,}\\\end{array}\right.
\]
\vspace{0.3cm}

\noindent then $Q=Q_1$ and $T=Q_0$. Hence by hypothesis, these sets are all nonempty and equal.

We will now prove that (5.2), (5.3), or (5.4) is satisfied by dividing the remainder of the argument into the three cases which are determined by the possible values of $k$.

\textbf{Case I}. Assume that $k=0$. We wish to verify the conclusion of (5.2). In this case $(2a_1, m)=1$, and so it follows from the exact form of the quadratic formula and the fact that $Q_0\not= \emptyset$ that $d/r^2\in Q(m)$. We next set $m_0=m$,
\[
\Sigma_i= \textrm{the set of all solutions of $x^2\equiv\ d/r^2$ mod  $m_i, i=0,1$,}\]
and let $u=\delta, v=m_1$ in Lemma 3.2 to conclude from that lemma and the equality $Q_0=Q_1$ that $\Sigma_0 \not= \emptyset \not= \Sigma_1$ and
\[
\Sigma _0=\big\{\sigma+jm_1: \sigma\in \Sigma_1,\ j\in \{0,\dots,\delta-1\}\big\}.\]
If we now let $s_i=$ the cardinality of $\Sigma_i,\ i=0, 1$, then $s_0\not=0\not=  s_1$ and
\begin{equation*}
\delta s_1=s_0. \tag{5.5}\]

For the next step in our argument, we will use the formula pointed out by Gauss that counts the number of solutions to the square-root problem. In order to state it, we let $u, v\in \textbf{Z}$ with $v>1$ and $u\in Q(v)$, consider the congruence
\begin{equation*}
x^2 \equiv\ u\ \textrm{mod}\ v,\tag{5.6}\]
and let
\[
\gamma=\textrm{the number of solutions of (5.6)}.\]
Suppose first that $v$ is a power $p^t$ of the prime $p$. It follows from Gauss' solution of the square-root problem that

$(a)$ If $u$ is not divisible by $p$ then $\gamma=1$ if $p=2$ and $t=1$, $\gamma=2$ if $p$ is odd or $p=2=t$, and $\gamma=4$ if $p=2, t>2$;

$(b)$ If $p$ divides $u$ and $p^t$ does not, let $2\mu=\mu_p(u)$ and set $u=u_1p^{2\mu}$. Then $\gamma=p^{\mu}, 2p^{\mu}$ or $4p^{\mu}$ if the number of solutions of $x^2 \equiv\ u_1\ \textrm{mod}\ p^{t-2\mu}$ is, respectively, 1, 2, or 4;

$(c)$ If $u$ is divisible by $p^t$ and $[\cdot]$ denotes the greatest integer function, then $\gamma=p^{[t/2]}$.

\noindent If $v$ is now an arbitrary modulus with prime factorization $p_1^{\alpha_1}\cdots p_s^{\alpha_s}$ and
\[
\gamma_i=\textrm{the number of solutions to $x^2 \equiv\ u$ mod $p_i^{\alpha_i}$},
\]
where $\gamma_i$ is calculated according to $(a), (b),$ or $(c), i=1,\dots,s$, then
\[
\gamma=\prod_{i=1}^s\ \gamma_i.\]

We next make three observations that will be of use to us momentarily:
\begin{equation*}
\textrm{if $p$ is an odd prime factor of $v$ which does not divide $u$ then $p$ is not a factor of $\gamma$;}\tag{5.7}\]
\begin{equation*}
\textrm{if $p$ is an odd prime factor of $v$ which divides $u$ then $\mu_p(\gamma)<\mu_p(v);$} \tag{5.8}\]
\begin{equation*}
\textrm{every odd prime factor of $\gamma$ is a factor of $v$.}\tag{5.9}\]

\vspace{0.1cm}
\noindent If $p$ is a prime factor of $\gamma$ then the multiplicity of $p$ in $\gamma$ will be called the \emph{counting multiplicity of p with respect to u and v}.

Consider now the prime factors of $m$. We divide them respectively into three sets $P_1, P_2,$ and  $P_3$: the prime factors of 
$\delta$ that are not factors of $m_1$, the common prime factors of $\delta$ and $m_1$, and the prime factors of $m_1$ that are not factors of $\delta$.

Assume that $\delta>1$. Let $p$ be a fixed prime factor of $\delta$. We will use the Gauss counting formula and equation (5.5) to analyze the multiplicity $\alpha$ of $p$ in $\delta$.

Begin by noting that $p$ is odd and a factor of the left-hand side of (5.5), hence also a factor of the right-hand side. We conclude by observation (5.7) that $d/r^2$ is divisible by $p$.

Suppose next that $p\in P_1$. Then $p$ is not a factor of $m_1$ and so $\alpha=\mu_p(m)$. Since $p$ is odd and not a factor of $m_1$, it follows from observation (5.9) that $p$ is not a factor of $s_1$. Hence $\alpha=\mu_p(\delta s_1)$. If $\mu(p)$ is the counting multiplicity of $p$ with respect to $d/r^2$ and $m$ then $\mu(p)=\mu_p(s_0)$. It follows that $\mu(p)=\alpha=\mu_p(m)$, and this contradicts observation (5.8). We conclude that $P_1$ is empty.

Suppose that $p \in P_2$. If $\beta=\mu_p(m_1)$ then $\alpha+\beta=\mu_p(m)$. If $\mu^{\prime}(p)$ denotes the counting multiplicity of $p$ with respect to $d/r^2$ and $m_1$, it follows from (5.5) that
\begin{equation*}
\alpha+\mu^{\prime}(p)=\mu(p). \tag{5.10}\]
If $p^{\beta}$ does not divide $d/r^2$, and if $2\mu=\mu_p(d/r^2)$, then $\mu^{\prime}(p)=\mu=\mu(p)$, which is not possible by (5.10). Hence $p^{\beta}$ divides $d/r^2$, and so $\mu^{\prime}(p)=[\beta/2]$. If $p^{\alpha+\beta}$ does not divide $d/r^2$ then $\mu(p)=\mu$. Now $2\mu$ does not exceed the largest even integer less than $\alpha+\beta$, hence
\[
\mu \leq \left[\frac{\alpha}{2} \right]+\left[\frac{\beta}{2} \right] .\]
Thus by (5.10),
\[
\alpha \leq \left[\frac{\alpha}{2} \right],\]
and no positive integer can satisfy this inequality. We conclude that $d/r^2$ is divisible by $p^{\alpha+\beta}$. Hence $\mu(p)=[(\alpha+\beta)/2]$, and so by (5.10),
\[
\alpha+\left[\frac{\beta}{2} \right]=\left[\frac{\alpha+\beta}{2} \right].\]
This equation implies that $\alpha=1$ and $\beta$ is odd, and so $p$ has even multiplicity in $m$. Hence $\delta$ is the product of distinct odd primes, every prime factor $p$ of $\delta$ has even multiplicity $m(p)$ in $m$, and $d/r^2$ is divisible by 
\[
\prod_{p\in P_2}\ p^{m(p)},\]
i.e., (5.2) is true.

\textbf{Case II}. We next suppose that $k=1$ and seek to verify the conclusion of (5.3). Take $u=2\delta,\ v=m_1$ in Lemma 3.2, let $S_i,\ \Sigma_i,\ i=0, 1$ be as defined in that lemma with this choice of $u$ and $v$, and thus conclude from the equality of $Q_0$ and $Q_1$ that $S_i,\ \Sigma_i,\ i=0, 1$ are nonempty and
\begin{equation*}
S_0=\big\{s+jm_1: s\in S_1,\ j\in \{0, 1,\dots,2\delta-1\} \big\}. \tag{5.11}\]

It follows that $d/r^2\in Q(8a_1m)$, and so $d/r^2\in Q(m)$, and from (5.11) and the fact that the cardinality of $S_i$ and $\Sigma_i$ are the same for $i=0, 1$, it also follows that
\begin{equation*}
\textrm{cardinality of $\Sigma_0=2\delta$(cardinality of $\Sigma_1$).} \tag{5.12}\]

Since $a_1$ and $m$ are odd and $(a_1, m)=1$, it is a consequence of Proposition 2.1 and the definition of $\Sigma_0$ that the elements of $\Sigma_0$ are obtained precisely as the simultaneous solutions $\sigma$ of
\[
\sigma \equiv \tau\ \textrm{mod}\ 8,\]
\[
\sigma \equiv \alpha\ \textrm{mod}\ a_1,\]
\[
\sigma \equiv \mu\ \textrm{mod}\ m,\]
where $\tau$ varies over all solutions of
\begin{equation*}
\tau^2 \equiv \frac{d}{r^2}\ \textrm{mod}\ 8 \tag{5.13}\]
that are pairwise incongruent mod 4, $\alpha$ varies independently over all solutions of
\begin{equation*}
\alpha^2 \equiv \frac{d}{r^2}\ \textrm{mod}\ a_1 \tag{5.14}\]
which also satisfy
\begin{equation*}
\alpha \equiv\ b_1 \textrm{mod}\ a_1, \tag{5.15}\]
and $\mu$ varies independently over all solutions of
\begin{equation*}
\mu^2 \equiv \frac{d}{r^2}\ \textrm{mod}\ m. \tag{5.16}\]
We note that if $d/r^2$ is even then there is exactly one such solution $\tau$, if $d/r^2$ is odd there are exactly two such solutions, and that $b_1$ always determines a solution of (5.14) and (5.15).

The same reasoning shows that the elements of $\Sigma_1$ consist precisely of the simultaneous solutions $\sigma$ of
\[
\sigma \equiv \tau\ \textrm{mod}\ 4,\]
\[
\sigma \equiv \alpha\ \textrm{mod}\ a_1,\]
\[
\sigma \equiv \mu\ \textrm{mod}\ m_1,\]
where $\tau$ varies over all solutions of
\begin{equation*}
\tau^2 \equiv \frac{d}{r^2}\ \textrm{mod}\ 4 \tag{5.17}\]
that are pairwise incongruent mod 2, of which there is only one such solution, $\alpha$ varies independently over all solutions of (5.14) and (5.15),
and $\mu$ varies independently over all solutions of
\begin{equation*}
\mu^2 \equiv \frac{d}{r^2}\ \textrm{mod}\ m_1. \tag{5.18}\]
It hence follows from (5.12) that if
\[
t=\textrm{cardinality of the set of all solutions of (5.13) that are pairwise incongruent mod 4,}\]
\[
s_0=\textrm{cardinality of the set of all solutions of (5.16),}\]
\[
s_1=\textrm{cardinality of the set of all solutions of (5.18),}\]
then
\begin{equation*}
2\delta s_1=ts_0.\tag{5.19}\]

Assume that $\delta>1$. Since $t$ is either 1 or 2, it follows that the analysis of $\delta$ that was carried out in the proof of (5.2) can also be done here, with (5.19) in place of (5.5), to show that $\delta$ and $d/r^2$ satisfy the conditions as specified for them in the conclusion of (5.2). But then $\delta s_1=s_0$, hence $t=2$, and so $d/r^2$ must be odd. Since $d/r^2\in Q(8)$, it hence follows that
\[
b_1^2-4a_1c_1=\frac{d}{r^2}\equiv 1\ \textrm{mod}\ 8,\]
and thus $b_1$ is odd and either $a_1$ or $c_1$ is even. If $\delta=1$ then $s_1=s_0$, hence $t=2$, and we conclude as before that $b_1$ is odd and either $a_1$ or $c_1$ is even in this case as well. We have verified (5.3).

\textbf{Case III}. Assume now that $k\geq 2$, and suppose first, by way of contradiction, that $\rho\geq k-2$. Let  $u=2^k\delta,\ v=m_1$ in Lemma 3.2 to conclude as before that if $\Sigma_i,\ i=0, 1$ are defined as in that lemma with this choice of $u$ and $v$, then these sets are nonempty,
\begin{equation*}
\textrm{cardinality of $\Sigma_0=2^k\delta$(cardinality of $\Sigma_1$)}, \tag{5.20}\]
and if $e_0=k+2,\ e_1=2$, and $m_0=m$, then the elements of $\Sigma_i$ are given by the simultaneous solutions of
\[
\sigma \equiv \tau\ \textrm{mod}\ 2^{e_i},\]
\[
\sigma \equiv \alpha\ \textrm{mod}\ a_1,\]
\[
\sigma \equiv \mu\ \textrm{mod}\ m_i,\]
where $\tau$ varies over all solutions of
\begin{equation*}
\tau^2 \equiv \frac{d}{r^2}\ \textrm{mod}\ 2^{e_i} \tag{5.21}\]
that are pairwise incongruent mod $2^{e_i-1}, \alpha$ varies independently over all solutions of (5.14) and (5.15),
and $\mu$ varies independently over all solutions of
\begin{equation*}
\mu^2 \equiv \frac{d}{r^2}\ \textrm{mod}\ m_i,\ i=0, 1. \tag{5.22}\]
If $s_0$ and $s_1$ are defined as in the proof of (5.3) and
\begin{quote}
$c_0=$ cardinality of the set of all solutions of (5.21) with $i=0$ that are pairwise incongruent mod $2^{k+1}$,
\end{quote}
then by (5.20),
\begin{equation*}
2^k\delta s_1=c_0s_0. \tag{5.23}\]

Our strategy here, as before, is to employ a counting argument which exploits (5.23). This requires the calculation of $c_0$. To that end, we first assert that 4 must divide $d/r^2$. In order to see that, let $x\in Q_1$ and deduce from the fact that $Q_0=Q_1$ that
\begin{equation*}
q_1(x+zm_1)\equiv 0\ \textrm{mod}\ 2^km,\ \forall z\in \textbf{Z}.\tag{5.24}\]
If we now use the fact that $q_1(x)\equiv 0\ \textrm{mod}\ 2^km (x\in Q_0$!) and take $z=2$ in (5.24), we obtain the congruence
\[
2a_1x+b_1+2a_1m_1\equiv 0\ \textrm{mod}\ 2,\]
i.e., $2a_1x+b_1$ is even. Since
\[
(2a_1x+b_1)^2\equiv \frac{d}{r^2}\ \textrm{mod}\ 4a_1m_1,\]
$d/r^2$ is hence divisible by 4.

Suppose now that $2^{k+2}$ does not divide $d/r^2$. If $2\mu=\mu_2(d/r^2)$, it is a straightforward consequence of Gauss' solution to the square-root problem that if $d/r^2=d_1\cdot 2^{2\mu}$, then the solutions of (5.21) for $i=0$ that are pairwise incongruent mod $2^{k+1}$ can be taken to be
\begin{equation*} 
\eta\cdot 2^{\mu}+s\cdot 2^{k+2-\mu},\ s\in \{0,\dots,2^{\mu-1}-1\},\tag{5.25}\]
where $\eta$ varies over all solutions of
\begin{equation*} 
\eta^2 \equiv d_1\ \textrm{mod}\ 2^{k+2-2\mu}.\tag{5.26}\]
Hence
\[
c_0=\varepsilon\cdot 2^{\mu-1},\]
where $\varepsilon=1, 2,$ or 4, depending on whether (5.26) has, respectively, 1, 2, or 4 solutions. Thus by (5.23),
\begin{equation*}
2^k\delta s_1=\varepsilon\cdot 2^{\mu-1} s_0.\tag{5.27}\]
If $\delta=1$ then $s_0=s_1$ and we obtain
\begin{equation*}
2^k=\varepsilon\cdot 2^{\mu-1}.\tag{5.28}\]
If $\delta>1$, we reason from (5.27) as in the proof of (5.3) to conclude that $\delta s_1=s_0$, and so we obtain (5.28) in this instance as well.

From (5.28) it follows that $k=\mu-1, \mu,$ or $\mu+1$. But each of these alternatives will occur if and only if $k+2-2\mu=1, k+2-2\mu=2$ or $k+2-2\mu\geq 3$, respectively, and so they can occur only if $\mu=0$, which is not possible.

We conclude that $d/r^2$ is divisible by $2^{k+2}$. Hence if $t$ is chosen so that $k+2=2t$ or $2t-1$, depending on the parity of $k+2$, then the solutions of (5.21) with $i=0$ which
are pairwise incongruent mod $2^{k+1}$ can be taken to be
\[
s\cdot 2^t,\ s\in \{0, 1,\dots,2^{k+1-t}-1\}.\]
Hence $c_0=2^{k+1-t}$ in this case, and so by (5.23),\[
2^k\delta s_1=2^{k+1-t} s_0.\]
By use of the same argument as before, this equation will be true only if $t=1$, i.e., $k=0$, contrary to hypothesis.

It follows that $\rho<k-2$. This situation now requires that we take $u=2^{\rho+2} \delta,\ v=2^{k-\rho-2} m_1$ in Lemma 3.2, define $\Sigma_0$ and $\Sigma_1$ as per that choice, note that $\Sigma_0\not= \emptyset \not= \Sigma_1$,
\begin{equation*}
\textrm{cardinality of $\Sigma_0=2^{\rho+2} \delta$(cardinality of $\Sigma_1$),}\tag{5.29}\]
\begin{equation*}
\textrm{for each $\sigma \in \Sigma_1$ and $j\in \{0, 1,\dots,2^{\rho+2} \delta-1\}$, there exists  $\sigma^{\prime} \in \Sigma_0$ such that}\tag{5.30}\] 
\[
\sigma^{\prime}\equiv \sigma+2^{k-\rho-1} a_1m_1j\ \textrm{mod}\ 2^{k+1} a_1m,\]
and that the elements of $\Sigma_0$ and $\Sigma_1$ are given by the simultaneous solutions of the same congruences as before via (5.14), (5.15), (5.21), and (5.22), with $e_0=k+2$ and $e_1=k-\rho$.

If $s_i$ is defined as before and
\begin{quote}
$c_i=$ cardinality of the set of all solutions of (5.21) that are pairwise incongruent mod $2^{e_i-1}, i=0, 1$,
\end{quote}
then we obtain via (5.29) that
\begin{equation*}
2^{\rho+2} \delta c_1s_1=c_0s_0.\tag{5.31}\]

We check that $d/r^2$ is still divisible by 4, and if we suppose that $2^{k-\rho}$ does not divide $d/r^2$, then straightforward modification of our previous reasoning show that if $2\mu=\mu_2(d/r^2)$ then
\[
c_0=2^{\mu+1},\ c_1=\varepsilon\cdot 2^{\mu-1},\]
where $\varepsilon=1, 2, $ or 4. We hence  conclude from (5.31) that $2^{\rho} \varepsilon=1$, i.e., $\rho=0$ and $\varepsilon=1$, in which case $k=2\mu+1$. It follows that
\begin{equation*}
\textrm{if $2^{k-\rho}$ does not divide $d/r^2$ then $\rho=0, k$ is odd, and $k-1=\mu_2(d/r^2)$.} \tag{5.32}\]

Suppose next that $d/r^2$ is divisible by $2^{k+2}$. Then $d/r^2$ is also divisible by $2^{k-\rho}$, and so if we choose $k+2$ (respectively, $k-\rho)=2s$ or $2s-1$ (respectively, $2t$ or $2t-1$), according to the relevant parities, we find that
\[
c_0=2^{k-s+1},\ c_1=2^{k-\rho-t-1},\]
hence from (5.31) it follows that $s=t$, obviously impossible. Thus 
\begin{equation*}
\textrm{$d/r^2$ is not divisible by $2^{k+2}$.}\tag{5.33}\]

We can now prove that $\rho=0, k$ is odd,and $k-1=\mu_2(d/r^2)$. In light of (5.32) this will be done by showing that $2^{k-\rho}$ does not divide $d/r^2$. In order  to do that, we observe first that from (5.30) it follows that
\begin{equation*}
\textrm{for each element $\tau$ of the set of solutions of (5.21) with $i=1$ and $e_1=k-\rho$ } \tag{5.34}
\]
\hspace{1.7cm} that are pairwise incongruent mod $2^{k-\rho-1}$ and $j\in \{0, 1,\dots,2^{\rho+2} \delta-1\}$, there 

\vspace{0.2cm}
\hspace{1.3cm} exits an element $\tau^{\prime}$ from the set of solutions of (5.21) with $i=0$ and $e_0=k+2$
\vspace{0.2cm}

\hspace{1.4cm} that are pairwise incongruent mod $2^{k+1}$ and $t\in \{0, 1,\dots,2^{\rho+1}-1\}$ such that
\begin{equation*}
\tau^{\prime}\equiv \tau+t\cdot 2^{k-\rho}+2^{k-\rho-1} a_1m_1j\ \textrm{mod}\ 2^{k+1}.\]

Suppose now that $d/r^2$ is divisible by $2^{k-\rho}$. Then if $k-\rho=2w$ or $2w-1$, the solutions of (5.21) with $i=1$ as in (5.34) can be taken to be
\begin{equation*}
s\cdot 2^w,\ s\in \{0,\dots,2^{k-\rho-w-1}-1\}.\tag{5.35}\]
By virtue of (5.33), if $2\mu=\mu_2(d/r^2)$ and $d/r^2=d_1\cdot 2^{2\mu}$, then the solutions of (5.21) with $i=0$ as in (5.34) can be taken as in (5.25) and (5.26).

Assume first that $\mu<k-\rho$. If we set $s=0$ in (5.35) and $j=2^{\rho+1}$ in (5.34), then we find $\eta$ as in (5.26) and integers $t$ and $u$ such that
\begin{equation*}
t\cdot 2^{k-\rho-\mu}+2^{k-\mu}a_1m_1\equiv \eta+u\cdot 2^{k+2-2\mu} \ \textrm{mod}\ 2^{k+1-\mu}.\tag{5.36}\]
Now it follows from (5.31) that $k-\mu\geq w-2$. Since $k-\rho>2,\ w$ must be at least 2, hence $k-\mu\geq 0$. But $k\not= \mu$ since $\rho$ is nonnegative. We thus conclude from (5.36) that $\eta$ is even, hence by(5.26) so is $d_1$, contradicting the fact that $2\mu=\mu_2(d/r^2)$.

We conclude that $k-\rho \leq \mu$. If $k-\rho \geq 4$ then we can take $s=1$ in (5.35) and $j=0$ in (5.34) to find integers $\eta, t$, and $u$ so that
\[
1+t\cdot 2^{k-\rho-w}\equiv 2^{\mu-w}\cdot \eta+u\cdot 2^{k+2-\mu-w}\ \textrm{mod}\ 2^{k+1-w}.
\]
Since $\mu\geq 2w-1,\ k+2-\mu>\mu,\ k-\rho-w\geq w-1,\ k+1-w\geq w$ and $w\geq 2$, this congruence yields another contradiction. Finally, if $k-\rho=3$, we must take $s=0$ in (5.35) and so if we choose $j=1$ in (5.34), we obtain integers $\eta, t$, and $u$ for which
\[
2t+a_1m_1\equiv 2^{\mu-2}\cdot \eta+u\cdot 2^{k-\mu}\ \textrm{mod}\ 2^{k-1}.\]
Because $w=2$, we have $\mu\geq 3$ and $k-\mu\geq \mu-1 \geq 2$, and since $a_1m_1$ is odd, this congruence also is impossible. It follows that $2^{k-\rho}$ does not divide $d/r^2$.

Because  $\rho=0, k$ is odd, and $k-1=\mu_2(d/r^2)$, it follows from (5.31) that $\delta s_1=s_0$ and so $\delta$ and $d/r^2$ satisfy the conditions specified in (5.2).

Finally, we deduce from the fact that $Q_0\not= \emptyset$ that  $d/r^2 \in Q(2^{k+2}a_1m)$, hence in particular, $d/r^2 \in Q(m) \cap Q(2^k)$. Now, as Gauss points out in [1, articles 102 and 103], the even integers in $Q(2^k)$ consist precisely of 0 and the integers $z$ which satisfy the following conditions: if $\mu=\mu_2(z)$ then either $\mu\geq k$ or $\mu$ is even, $0<\mu<k$, and $z/2^{\mu} \equiv 1$ mod 8. As $k$ is odd, $k-1=\mu_2(d/r^2)$, and $d/r^2 \in Q(2^k)$, it hence follows that $d/(r^2\cdot 2^{k-1})\equiv 1$ mod 8. We have verified (5.4). $\hspace{13cm}\textrm{QED}$

We now state and prove the converse of Lemma 5.1, after retaining the notation as specified in the statement of that lemma.
\begin{lem}
If condition $(5.1)$ holds and either $k=0$ and the conclusion of $(5.2)$ holds, or $k=1$ and the conclusion of $(5.3)$ holds, or $k\geq 2$ and the conclusion of $(5.4)$ holds, then $\emptyset\not= Q=T$.
\end{lem}

\emph{Proof}. Suppose (5.1) is true. If $Q_0$ and $Q_1$ are defined as before then $Q=Q_1$ and $T=Q_0$, and so under each of the hypotheses in Lemma 5.2, we must prove that
\begin{equation*}
\emptyset\not= Q_1=Q_0.\tag{5.37}\]
As in the proof of Lemma 5.1, we divide the reasoning into the cases which are determined by the possible values of $k$.

\textbf{Case IV}. Assume to begin with that $k=0$ and the conclusion of (5.2) is true. Because $(2a_1, m)=1$ and $d/r^2\in Q(m)$, we conclude from the exact form of the quadratic formula that $Q_0\not= \emptyset\not= Q_1$.

If $\delta=1$ then $m=m_1$, and so (5.37) is trivially true. Hence assume that $\delta>1$; then $\delta$ also satisfies the conditions as specified in (5.2). Letting $\Sigma_0$ and $\Sigma_1$ be defined as before in this case, we have that $\Sigma_0 \not= \emptyset \not= \Sigma_1$ , hence we must prove, by virtue of Lemma 3.2, that
\begin{equation*}
\Sigma _0=\big\{\sigma+jm_1: \sigma\in \Sigma_1,\ j\in \{0,\dots,\delta-1\}\big\}.\tag{5.38}\]
Since $\Sigma_0$ is clearly contained in the set on the right-hand side of (5.38), we need only verify the reverse inclusion.

Let $p_1^{\beta_1}\cdots p_t^{\beta_t}$  be the prime factorization of $m_1$. It follows from the conditions satisfied by $\delta$ that $m$ and $m_1$ have the same prime factors, and if  $p_1^{\alpha_1}\cdots p_t^{\alpha_t}$ is the prime factorization of $m$, then whenever $p_i$ is a common prime factor of $\delta$ and $m_1$, we have that $\alpha_i$ is even, $\beta_i=\alpha_i-1$ and $p_i^{\alpha_i}$ divides $d/r^2$, and whenever $p_i$ is a factor of $m_1$ that is not a factor of $\delta$, then $\alpha_i=\beta_i$. As in the proof of Lemma 5.1, let $P_2$ and $P_3$ denote, respectively, the set of common prime factors of $\delta$ and $m_1$ and the set of prime factors of $m_1$ which are not factors of $\delta$.

Let $\sigma\in \Sigma_1$ and $j\in \{0, 1,\dots, \delta-1\}$; we will find $\sigma^{\prime}\in \Sigma_0$ such that
\begin{equation*}
\sigma^{\prime}\equiv \sigma+jm_1\ \textrm{mod}\ m.\tag{5.39}\]
In order to do that, we first find a solution $x_i$ of $x^2\equiv d/r^2\ \textrm{mod}\ (p_i^{\beta_i})$ such that
\begin{equation*}
\sigma \equiv x_i\ \textrm{mod}\ p_i^{\beta_i},\ i=1,\dots,t.\tag{5.40}\]
Next, for each prime $p_i\in P_2$, we find $q_i\in \textbf{Z}$ such that
\begin{equation*}
\sigma \equiv x_i+q_ip_i^{\beta_i}\ \textrm{mod}\ p_i^{\alpha_i}.\tag{5.41}\]
We now claim that
\begin{equation*}
\textrm{for each $p_i\in P_2$, there exists a solution $x_i^{\prime}$ of $x^2\equiv d/r^2\ \textrm{mod}\ p_i^{\alpha_i}$ such that}\tag{5.42}\]
\[
x_i^{\prime}\equiv x_i+q_ip_i^{\beta_i}+jm_1\ \textrm{mod}\ p_i^{\alpha_i}.\]
If (5.42) is true then we find $\sigma^{\prime}\in \Sigma_0$ such that
\begin{equation*}
\sigma^{\prime}\equiv x_i^{\prime}\ \textrm{mod}\ p_i^{\alpha_i},\ \textrm{if}\ p_i\in P_2,\tag{5.43}\]
\begin{equation*}
\sigma^{\prime}\equiv x_i\ \textrm{mod}\ p_i^{\alpha_i},\ \textrm{if}\ p_i\in P_3.\tag{5.44}\]
After observing that $m_1$ is divisible by $p_i^{\alpha_i}$ whenever $p_i\in P_3$, it follows from (5.40)-(5.44) that
\[
\sigma^{\prime}\equiv \sigma+jm_1\ \textrm{mod}\ p_i^{\alpha_i},\ i=1,\dots,t,\]
and this yields (5.39).

In order to establish (5.42), we fix $p_i=p\in P_2$, set $q=q_i$ and let $\alpha_i=2s,\ \beta_i=2s-1$. Since $d/r^2$ is divisible by $p^{2s}$, the solutions of $x^2\equiv d/r^2\ \textrm{mod}\ (p^{2s})$ and $x^2\equiv d/r^2\ \textrm{mod}\ (p^{2s-1})$ are given, respectively, by
\begin{equation*}
\big\{ip^s: i\in \{0, 1, \dots,p^s-1\}\big\},\tag{5.45}\]
\begin{equation*}
\big\{ip^s: i\in \{0, \dots,p^{s-1}-1\}\big\}.\tag{5.46}\]

Let $i\in \{0, \dots,p^{s-1}-1\}$. Then in view of (5.45) and (5.46), (5.42) will be true if we can find $v\in \{0, 1, \dots,p^s-1\}$ such that
\begin{equation*}
ip^s+qp^{2s-1}+jm_1\equiv vp^s\ \textrm{mod}\ p^{2s}.\tag{5.47}\]
But that can be done by first observing that $m_1$ is divisible by $p^s$, and so there is a $u\in \textbf{Z}$ such that
\begin{equation*}
m_1\equiv up^s\ \textrm{mod}\ p^{2s}.\tag{5.48}\]
Now simply choose $v\in \{0, 1, \dots,p^s-1\}$ such that
\[i+qp^{s-1}+ju\equiv v\ \textrm{mod}\ p^{s},\]
multiply this congruence by $p^s$, and substitute (5.48) into the congruence that results to obtain (5.47). This verifies (5.38).

\textbf{Case V}. Suppose next that $k=1$ and the conclusions of (5.3) are true. Since $b_1$ is odd and either $a_1$ or $c_1$ is even, it follows that $d/r^2\equiv 1$ mod 8, hence $d/r^2\in Q(8)$. Since by hypothesis we also have that $d/r^2\in Q(m)$, and $b_1$ determines a solution of (5.14) and (5.15), it is a consequence of the recipe for the construction of the elements of $\Sigma_0$ and $\Sigma_1$ for this case that  $\Sigma_0\not= \emptyset \not=\Sigma_1$. In order to verify (5.37), we must, as per Lemma 3.2, show that for each $\sigma \in \Sigma_1$ and $j\in \{0, 1,\dots,2\delta-1\}$, there exists $\sigma^{\prime} \in \Sigma_0$ such that
\[
\sigma^{\prime}\equiv \sigma+2a_1m_1j \ \textrm{mod}\ 4a_1m,\]
and this will hold if we in turn prove that
\begin{equation*}
\textrm{for each element $\tau$ from the set of solutions of (5.17) that are pairwise incon-} \tag{5.49}\] 

\hspace{1.3cm} grunt mod 2 and each $j\in \{0, 1,\dots,2\delta-1\}$, there exists an element $\tau^{\prime}$ from 

\vspace{0.2cm}
\hspace{1.3cm} the set of solutions of (5.13) that are pairwise incongruent mod 4 such that
\[
\tau^{\prime}\equiv \tau+2a_1m_1j \ \textrm{mod}\ 4,\]
and
\begin{equation*}
\textrm{for each solution $\mu$ of (5.18) and  $j\in \{0, 1,\dots,2\delta-1\}$, there exists a solution}\tag{5.50}\]

\hspace{1.4cm}$\mu^{\prime}$ of (5.16) such that
\[
\mu^{\prime}\equiv \mu+2a_1m_1j \ \textrm{mod}\ m.\]

It follows from the hypothesis on $\delta$ and our previous reasoning that (5.50) is valid. In order to verify (5.49), we first observe that $d/r^2$ is odd, hence in (5.49) $\tau$ is either 1 or 3 and $\tau^{\prime}$ is either 1 or 3, 1 or 7, 3 or 5, or 5 or 7. Thus for any allowable $\tau$ and $j$,
\[
\tau+2a_1m_1j\equiv 1\ \textrm{or}\ 3 \ \textrm{mod}\ 4,\]
and so there is an appropriate $\tau^{\prime}$ which makes (5.49) true.

\textbf{Case VI}. Suppose finally that $k\geq 2$ and the conclusion of (5.4) is true. Because $k$ is odd, $k-1=\mu_2(d/r^2)$, and $d/(r^2\cdot 2^{k-1})\equiv 1$ mod 8, it follows that $d/r^2\in Q(2^{k+2})$. This together with the assumption $d/r^2\in Q(m)$ implies that  $\Sigma_0\not= \emptyset \not=\Sigma_1$ in this case. Hence we must prove that for each $\sigma \in \Sigma_1$ and $j\in \{0, 1,\dots,4\delta-1\}$, there exists $\sigma^{\prime} \in \Sigma_0$ such that
\[
\sigma^{\prime}\equiv \sigma+2^{k-1}a_1m_1j \ \textrm{mod}\ 2^{k+1}a_1m,\]
and this in turn will be so if (5.50) holds with $4\delta$ and $2^{k-1}a_1m_1j$ in place of $2\delta$ and $2a_1m_1j$, respectively, and if
\begin{equation*}
\textrm{for each $\varepsilon\in \{0, 1\}$, for each $j\in \{0, 1,\dots,4\delta-1\}$, and for each element $\tau$ from }\tag{5.51}\]

\hspace{1.05cm} the set of solutions of (5.21) with $i=1$ and $e_1=k$ that are pairwise incongru- 

\vspace{0.2cm}
\hspace{1.05cm} net mod $2^{k-1}$, there exists an element $\tau^{\prime}$ from the set of solutions of (5.21) with

\vspace{0.2cm}
\hspace{1.05cm} $i=0$ and $e_0=k+2$ that are pairwise incongruent mod $2^{k+1}$ such that
\[
\tau^{\prime}\equiv \tau+\varepsilon \cdot 2^k+2^{k-1}a_1m_1j \ \textrm{mod}\ 2^{k+1}.\]
But (5.50) as modified holds by the same reasoning as before, so we need only verify (5.51). To that end, let $2\mu=\mu_2(d/r^2), d_1=d/(2^{2\mu}r^2)$, and so $k=2\mu+1.$ Verification of (5.51) requires showing that for each $i\in \{0,\dots,2^{\mu-1}-1\}, \varepsilon \in \{0, 1\},$ and $j\in \{0,1,\dots,4\delta-1\}$, there is an $s\in \{0,\dots,2^{\mu-1}-1\}$ and a solution $\eta$ of $\eta^2\equiv d_1$ mod 8 such that
\begin{equation*}
2^{\mu}+i\cdot 2^{\mu+1}+\varepsilon\cdot 2^{2\mu+1}+2^{2\mu}a_1m_1j\equiv \eta\cdot 2^{\mu}+s\cdot 2^{\mu+3}\ \textrm{mod}\ 2^{2\mu+2}.\tag{5.52}\]
Because $d_1\equiv$ 1 mod 8, $\eta$ can be either 1, 3, 5, or 7, hence this congruence will be satisfied for $i, \varepsilon, j, s$, and $\eta$ as specified if there exist an $s$ as specified and $\eta^{\prime}\in \{0, 1, 2, 3\}$ such that
\[
i+\varepsilon\cdot 2^{\mu}+2^{\mu-1}a_1m_1j\equiv \eta^{\prime}+4s\ \textrm{mod}\ 2^{\mu+1}.\]
Observe now that as $\eta^{\prime}$ and $s$ vary independently over all elements of 
\[
\{0, 1, 2, 3\}\ \textrm{and}\ \{0,\dots,2^{\mu-1}-1\},\] 
respectively, $\eta^{\prime}+4s$ varies over all elements of $ \{0, 1,\dots,2^{\mu+1}-1\}$, and this last set is a complete set of residues mod $2^{\mu+1}$. If $i, \varepsilon$, and $j$ are chosen as specified it thus follows that an appropriate $\eta$ and $s$ can be found so that (5.52) is true. Hence (5.37) is also true. QED
\section{The Main Theorem, Corollaries, and Examples}

Lemmas 3.6, 4.1, 5.1, and 5.2 now supply a proof of the following theorem, the principal result of this paper.
\begin{thm}
Let $a, b, c, n\in \textbf{Z}$, with $n\geq 2$ and a not divisible by n. If $r=(a, n),\ k=$ multiplicity of $2$ in $n/r,\ m=n/(2^kr),\ \delta=(m, r)$, and $d=b^2-4ac$ then $\textnormal{IQF}$ is valid for $ax^2+bx+c\equiv 0\ \textnormal{mod}\ n$ if and only if either

$(a)$ d is a quadratic non-residue of n, or

$(b)$ d is a quadratic residue of n and exactly one of the following mutually exclusive conditions holds:

$(i)$ there exists a prime factor p of $2a$ such that if $\beta=$ multiplicity of p in $2a$ and
\[
\alpha=\left\{\begin{array}{cc} \textrm{multiplicity of $p$ in $4am$, if $k=0, 1,$ or $2$,}\\
\textrm{multiplicity of $p$ in $2^kam$, if $k\geq$ $3$},\\\end{array}\right.
\]
then $1<\beta<\alpha,$ b is divisible by p, d is divisible by $p^2$, and $(p^{\alpha}, p^{\beta})$ forms a $(b, d)$-obstruction;

$(ii)\ k=0,$ b and c are divisible by r, $d/r^2$ is a quadratic residue of m and either $\delta=1$ or $\delta$ is the product of distinct odd primes $p_1,\dots, p_t,$ each prime $p_i$ has even multiplicity $m_i$ in $m$, and $d/r^2$ is divisible  by the product $p_1^{m_1}\cdots p_t^{m_t}$;

$(iii)\ k=1,$ b and c are divisible by r, $b/r$ is odd, either $a/r$ or $c/r$ is even, and $d/r^2$ and $\delta$ satisfy the conditions specified for them in $(b)(ii)$;

$(iv)\ k \geq 3$, r and k are odd, b and c are divisible by r, $k-1$ is the multiplicity of $2$ in $d/r^2,\ d/(r^2\cdot 2^{k-1})\equiv1\ \textnormal{mod}\ 8$, and $d/r^2$ and $\delta$ satisfy the conditions specified for them in $(b)(ii)$.
\end{thm}

\textbf{Remark}. The condition ``$d/r^2$ is a quadratic residue of $m$" in Theorem 6.1$(b)(ii)$-$(iv)$ may be replaced there by the condition ``$d/r^2$ is a quadratic residue of $m/\delta$".

The following corollaries of Theorem 6.1 give necessary and sufficient conditions for the validity of IQF in the interesting special case of a prime-power modulus. We note incidentally that if $p$ is an odd prime and $i\in \textbf{Z}^{+}$, then $(2a, p^i)=1$ if and only if $(a, p)=1$, hence we may suppose that $(a, p)>1$ in this case. We also maintain the notation used in  Theorem 6.1.
\begin{cor}
Let p be an odd prime, $i\in \textbf{Z},\ i\geq 2$. If $(a, p^i)=p^l,\ 1\leq l<i$, then $\textnormal{IQF}$ is valid for 
\[
ax^2+bx+c\equiv 0\ \textnormal{mod }p^i \]
if and only if either

$(a)$ d is a quadratic non-residue of $p^i$ or

$(b)$ d is a quadratic residue of $p^i$ and exactly one of the following mutually exclusive conditions holds:

$(i)\ l>1,$ b is divisible by p, and $(p^i, p^l)$ forms a $(b, d)$-obstruction;

$(ii)\ l=1,$ i is odd, b and c are divisible by p, and d is divisible by $p^{i+1}$.
\end{cor}

\emph{Proof}. We have that $m=p^{i-l}$ and $k=0$ in the hypotheses of Theorem 6.1. Thus IQF is valid for $ax^2+bx+c\equiv 0$ \textnormal{mod }$p^i$ if and only if $(a), (b)(i),$ or $(b)(ii)$ of that theorem holds.

Let $q$ be a prime factor of $2a$ and let $\alpha=\mu_q(4am)=\mu_q(4ap^{i-l}),\ \beta=\mu_q(2a)$ and $\mu=\mu_2(a)$. If $q=2$ then $\alpha=\mu+2,\ \beta=\mu+1$, and if $q=p$ then $l=\mu_p(a)$, and so $\alpha=i,\ \beta=l$. If $p\not= q \not= 2$ then $\alpha=\mu_q(a/p^l)=\beta$. It follows that condition $(b)(i)$ of Theorem 6.1 can hold only if the prime there is either 2 or $p$. We will prove that it cannot be 2.

Suppose that it is. Then, in particular, $(2^{\mu+2}, 2^{\mu+1})$ forms a $(b, d)$-obstruction, i.e., either $(b)(iii)$ or $(b)(iv)$ of Proposition 3.3 must hold for this pair.

Assume that $(b)(iii)$ of Proposition 3.3 holds. Then $2^{\mu+2}$ does not divide $d$, and if $2\nu =\mu_2(d),\ d_1=d/2^{2\nu}$, and $\Sigma$ is the set of all solutions of $x^2\equiv d_1$ mod $2^{\mu+2-2\nu}$, then
\begin{equation*}
b\not \equiv \sigma\cdot2^{\nu}+i\cdot 2^{\mu+2-\nu}\ \textrm{mod}\ 2^{\mu+1},\ \forall\ \sigma\in \Sigma,\ \forall i\in \{0, 1\dots, 2^{\nu}-1\}.\tag{6.1}\]

As $\mu+2=\mu_2(4a),\ b^2\equiv d$ mod $4a$, and $2^{\nu}<\mu+2$, it follows that $2\nu=\mu_2(b^2)$, and so $\nu=\mu_2(b)$. Letting $b_1=b/2^{\nu}$, we conclude that
\[
b_1^2\equiv d_1\ \textrm{mod}\ 2^{\mu+2-2\nu}.\]
Hence there exists $\sigma\in \Sigma$ such that
\[
b^2\equiv \sigma\cdot 2^{\nu}\ \textrm{mod}\ 2^{\mu+2-\nu}.\]
But then for some $i\in \{0, 1\dots, 2^{\nu}-1\}$,
\[
b\equiv \sigma\cdot 2^{\nu}+i\cdot2^{\mu+2-\nu}\ \textrm{mod}\ 2^{\mu+2},\]
and this contradicts (6.1).

We conclude that $(b)(iv)$ of Proposition 3.3 must hold, i.e., $2^{\mu+2}$ divides $d$, and if $\mu+2=2s$ or $2s-1$ then
\begin{equation*}
b\not \equiv  i\cdot 2^s\ \textrm{mod}\ 2^{\mu+1},\ \forall i\in \{0, 1,\dots,2^{\mu+2-s}-1\}.\tag{6.2}\]
But $2^{\mu+2}$ also divides $b^2$, hence $\mu_2(b)\geq s$, and so we can find $ i\in \{0, 1,\dots,2^{\mu+2-s}-1\}$ such that\[
b \equiv  i\cdot 2^s\ \textrm{mod}\ 2^{\mu+2},\]
which contradicts (6.2). It now follows that either $(a)$ or $(b)(i)$ of Theorem 6.1 holds if and only if $(a)$ or $(b)(i)$ of Corollary 6.2 holds.

We determine next when $(b)(ii)$ of Theorem 6.1 is valid. We have $r=p^l$, so $\delta=(m, r)=p^{\min \{l, i-l\}}>1$. Hence $\delta$ is a product of distinct prime factors and every prime factor of $\delta$ has even multiplicity in $m=p^{i-l}$ if and only if $l=1$ and $i$ is odd. But if $l=1$ then $r=p$, and so the remaining requirements in $(b)(ii)$ of Theorem 6.1 will hold if and only if $b$ and $c$ are divisible by $p$ and $d$ is divisible by $p^{i+1}$. Thus $(b)(ii)$ of Theorem 6.1 is valid if and only if $(b)(ii)$ of Corollary 6.2 is also. $\hspace{10.1cm} \textrm{QED}$

\begin{cor}
If $i\in \textbf{Z}^{+}$ and $(a, 2^i)=2^l,\ l<i$ then $\textnormal{IQF}$ is valid for 
\[
ax^2+bx+c\equiv 0\ \textnormal{mod } 2^i\] 
if and only if either

$(a)$ d is a quadratic non-residue of $2^i$ or

$(b)$ d is a quadratic residue of $2^i$ and exactly one of the following mutually exclusive conditions holds:

$(i)\ l>0,\ i\geq l+3,\ b$ is even, and $(2^i, 2^{l+1})$ forms a $(b, d)$ obstruction;

$(ii)\ i=l+1$, b and c are divisible by $2^l,\ b/2^l$ is odd, and $c/2^l$ is even;

$(iii)$ a is odd, i is odd and at least $3$, and $i-1$ is the multiplicity of $2$ in d.
\end{cor}

\emph{Proof}. We have $r=2^l,\ m=\delta=1$, and $k=i-l\geq 1$ in the hypotheses of Theorem 6.1, and so IQF is valid for $ax^2+bx+c\equiv 0$ \textnormal{mod }$2^i$ if and only if $(a),\ (b)(i),\ (b)(iii)$, or $(b)(iv)$ of that theorem holds. 

Suppose that $l=0$, i.e., $a$ is odd. If $t\in \textbf{Z}^{+},\ p$ is a prime factor of $2a,\ \alpha=\mu_p(2^ta)$, and $\beta=\mu_p(2a)$, then either $\beta=1$ (if $p=2$) or $\alpha=\beta$ (if $p$ is odd). Hence $(b)(i)$ of Theorem 6.1 cannot hold in this case.

If $t\in \textbf{Z}^{+}$ and $p$ is an odd prime factor of $2a$ then $\mu_p(2a)=\mu_p(2^ta)$, and so  $(b)(i)$ of Theorem 6.1 will be valid only if $l>0$ and the prime there is 2. Since $l=\mu_2(a)$, we have in this case that
\[
1<l+1=\mu_p(2a)<l+t=\mu_2(2^ta),\ t\geq 2.\]
Consequently, if $i=l+1$ or $l+2$ then $\alpha=l+2,\ \beta=l+1$ in $(b)(i)$ of Theorem 6.1, and so this condition can hold only if  $(2^{l+2}, 2^{l+1})$ forms a $(b, d)$-obstruction, which is impossible, as we showed in the proof of Corollary 6.2. We conclude that  $(b)(i)$ of Theorem 6.1 is equivalent to condition $(b)(i)$ of Corollary 6.3, and we clearly have 
$(b)(iii)$ of Theorem 6.1 and $(b)(ii)$ of Corollary 6.3 equivalent.

If $i\geq l+2$ then $(b)(iv)$ of Theorem 6.1 is true if and only if $l=0,\ i$ is odd, $i-1=\mu_2(d)$, and $d/2^{i-1}\equiv 1$ mod 8, and this is equivalent to $(b)(iii)$ of Corollary 6.3.$\hspace{2.6cm} \textrm{QED}$

We close our discussion with the following table, which lists some simple examples of congruences $q(x)\equiv 0$ mod $n$ for which IQF is valid, and shows that none of the conditions stated in Theorem 6.1 or Corollary 6.2 or 6.3 can be deleted.

\vspace{0.4cm}
\begin{center}
\textbf{Table 1. Examples of IQF}
\end{center}

\vspace{0.4cm}
\begin{center}
\begin{tabular}{| l | r |  p{9cm} |}
\hline
$q(x)$ & $n$ & Justification of IQF\\
\hline
$3x^2+1$ & 9 & Theorem 6.1$(a)$, Corollary 6.2$(a)$\\
\hline
$x^2+x+1$ & 8 & Theorem 6.1$(a)$, Corollary 6.3$(a)$\\
\hline
$18x^2+18x+1$ & 27 & Theorem 6.1$(b)(i)$, Corollary 6.2$(b)(i)$ ($(b)(iii)$ of Proposition 3.3 satisfied)\\
\hline
$9x^2+3x+1$ & 27 & Theorem 6.1$(b)(i)$, Corollary 6.2$(b)(i)$ ($(b)(iv)$ of Proposition 3.3 satisfied)\\
\hline
$8x^2+2x+1$ & 64 & Theorem 6.1$(b)(i)$, Corollary 6.3$(b)(i)$ ($(b)(iii)$ of Proposition 3.3 satisfied)\\
\hline
$3x^2+6x+3$ & 27 & Theorem 6.1$(b)(ii)$, Corollary 6.2$(b)(ii)$\\
\hline
$x^2+x$ & 2 &  Theorem 6.1$(b)(iii)$, Corollary 6.3$(b)(ii)$\\
\hline
$x^2+2x$ & 8 & Theorem 6.1$(b)(iv)$, Corollary 6.3$(b)(iii)$\\
\hline

\end{tabular}
\end{center}

\vspace{1cm}
\begin{center}
R\textsc{eference}
\end{center}

\vspace{0.2cm}
[1] C. F. Gauss, \emph{Disquisitiones Arithmeticae}, 2nd edition, Dietrich, G$\ddot{\textrm{o}}$ttingen, 1870; 

\hspace{0.42cm} \emph{Werke}, Band I, Georg Olms Verlag, Hildescheim, 1973: English translation by A. A. 

\hspace{0.44cm} Clarke, Springer-Verlag, New York-Berlin-Heidelberg, 1986.

\end{document}